\documentclass{elsart3-1}

\usepackage{amsmath}
\usepackage{amssymb}
\usepackage{color}

\usepackage[english,francais]{babel}

\newtheorem{theorem}{Theorem}[section]
\newtheorem{lemma}[theorem]{Lemma}
\newtheorem{e-proposition}[theorem]{Proposition}

\newtheorem{e-definition}[theorem]{Definition\rm}

\newtheorem{theoreme}{Th\'eor\`eme}[section]

\newtheorem{proposition}[theoreme]{Proposition}

\renewcommand{\theequation}{\arabic{equation}}
\newcommand{\ai}[1]{#1_{\infty}}
\DeclareMathOperator{\Fi}{F_{\infty}}
\DeclareMathOperator{\supp}{supp}

\def\og{\leavevmode\raise.3ex\hbox{$\scriptscriptstyle\langle\!\langle$~}}
\def\fg{\leavevmode\raise.3ex\hbox{~$\!\scriptscriptstyle\,\rangle\!\rangle$}}

\journal{the Acad\'emie des sciences}
\begin{document}
\centerline{}
\begin{frontmatter}
\selectlanguage{english}
\title{On a price formation free boundary model by Lasry \& Lions: The Neumann problem}
\selectlanguage{english}
\author[authorlabel1]{Luis A. Caffarelli}
\ead{caffarel@math.utexas.edu}
\author[authorlabel2,authorlabel3]{Peter A. Markowich}
\ead{P.A.Markowich@damtp.cam.ac.uk}
\author[authorlabel3]{Marie-Therese Wolfram}
\ead{marie-therese.wolfram@univie.ac.at}
\address[authorlabel1]{Department of Mathematics, Institute for Computational Engineering and Sciences, University of Texas at Austin, USA}
\address[authorlabel2]{DAMTP, University of Cambridge, Cambridge CB3 0WA, UK}
\address[authorlabel3]{Faculty of Mathematics, University of Vienna, 1090 Vienna, Austria}
\medskip
\begin{center}
{\small Received *****; accepted after revision +++++\\
Presented by }
\end{center}
\begin{abstract}
\selectlanguage{english}
We discuss local and global existence and uniqueness for the price formation free boundary model with homogeneous Neumann boundary conditions introduced by Lasry \& Lions in 2007. The results are based on a transformation of the problem to the heat equation with nonstandard boundary conditions. The free boundary becomes the zero level set of the solution of the heat equation. The transformation allows us to construct an explicit solution and discuss the behavior of the free boundary. Global existence can be verified under certain conditions on the free boundary and examples of non-existence are given.
\vskip 0.5\baselineskip
\selectlanguage{francais}
\noindent{\bf R\'esum\'e} \vskip 0.5\baselineskip \noindent
Nous discutons l'existence locale et globale, ainsi que l'unicit\'e des solutions pour le mod\`ele de formation des prix \`a fronti\`ere libre avec des conditions aux bords de Neumann homog\`enes introduit par Lasry \& Lions en 2007. Nos r\'esultats sont bas\'es sur une transformation de ce probl\`eme en une \'equation de la chaleur avec des conditions aux bords non standard. La fronti\`ere libre devient la ligne de niveau z\'ero de la solution de l'\'equation de la chaleur. Cette transformation nous permet de construire une solution explicite et de discuter le comportement de la fronti\`ere libre. L'existence globale peut \^etre v\'erifi\'ee sous certaines conditions sur la fronti\`ere libre, et nous donnons des exemples de non-existence.
\end{abstract}
\end{frontmatter}
\selectlanguage{english}
\section{Introduction}
\label{sec:intro}
In this paper we discuss the mathematical analysis of a price formation mean field model, as stated in \cite{Lasry2007}. The equation models the trading of an economic good between a large group of buyers and a large group of vendors. It consists of a non-linear parabolic free boundary evolution equation which describes the evolution of the densities of buyers and vendors and also determines the price of the good. Here we consider the problem on the real interval $(-L,L)$, where $L$ denotes the maximum  and $-L$ the minimum price. The model is given by the free boundary value problem
\begin{subequations}\label{eq:pfm}
 \begin{align}
 f_t - f_{xx} &= \lambda(t)(\delta(x-p(t)+a) - \delta(x-p(t)-a)),\, x\in(-L,L),t>0,\label{eq:pfm1}\\
\lambda(t) &= -f_x(p(t),t),\; f(p(t),t)=0,\label{eq:pfm2}\\ 
 f(x,0) &= f_I,\, p(0) = p_0,\;\text{for some $p_0$ in $(-L+a,L-a)$},\label{eq:initial}\\
f_x(\pm L,t) &= 0,~t>0, 
\end{align}
\end{subequations}
with $ 0 < a < L$ and compatibility conditions at time $t=0$: $ f_I(p_0) = 0\text{ and }f_I(x) > 0\text{ for }x < p_0 \text{ and }f_I(x) < 0\text{ for }x > p_0$. Throughout this paper we assume that $f_I$ is in $L^2(-L,L)$ and denote the $L^2$ norm on $(-L,L)$ by $\lVert \cdot \rVert$. System \eqref{eq:pfm} has been studied in a number of papers, see \cite{Chayes2009,Markowich2009,Caffarelli2011}. Here we present global (non-)existence results of a solution of  \eqref{eq:pfm} on the bounded interval $(-L,L)$.
\section{Analysis of the Neumann problem - Transformation to the Heat Equation}
In the sequel we denote the positive and negative part of a function $f$ defined for $x \in (-L,L)$ by $f^+:= \max(f,0)$, $f^-:= \max(-f,0)$. We extend $f^+$,$f^-$ by the value $0$ outside the interval $(-L,L)$.
\begin{proposition}
Let $f$ be a solution of the modified price formation problem consisting of equation \eqref{eq:pfm1} with conditions \eqref{eq:pfm2}, initial datum \eqref{eq:initial} and modified boundary conditions (bc):
\begin{align}\label{eq:newbc}
~~f_x(-L,t) = 
\begin{cases}
0 \quad & p(t) > -L+a\\
f_x(-L+a,t) \quad & p(t) \leq -L+a.
\end{cases}
f_x(L,t) = 
\begin{cases}
0 \quad & p(t) < L-a\\
f_x(L-a,t) & p(t) \geq L-a.
\end{cases}
\end{align}
Then the function
\begin{align}\label{eq:trans}
F(x,t) = 
\begin{cases}
\phantom{-}\sum\nolimits^{\infty}_{n=0} f^+(x+na,t), x < p(t)\\
-\sum\nolimits ^{\infty}_{n=0} f^-(x-na,t), x > p(t),
\end{cases}
\end{align}
is a solution of the following BVP for the heat equation
\begin{subequations}\label{eq:heat}
\begin{align}
F_t &=  F_{xx}, ~\text{ for all } x \in (-L,L), t>0,\\
F_x(\pm L,t) &= F_x(\pm L\mp a,t), ~t>0\\
F(x,t=0) &= F_I(x), ~\text{ for all } x \in (-L,L),
\end{align}
\end{subequations}
where $F_I$ is constructed from $f_I$ according to \eqref{eq:trans}. Conversely, if $F$ is a solution of \eqref{eq:heat} then $f(x,t) = F(x,t)-F^+(x+a,t)+F^-(x-a,t)$ satisfies \eqref{eq:pfm1}, \eqref{eq:pfm2}, \eqref{eq:initial}, \eqref{eq:newbc}.
\end{proposition}
The construction of the function $F$ is motivated by the fact that $f^+$ ($f^-$) has jump-discontinuities at $x = p(t)-a$ ($x = p(t)+a$) and $x = p(t)$ of equal magnitude but opposite signs. The summation procedure then moves the jumps in the derivatives of $f$ outwards. Note that the sum in \eqref{eq:trans} actually consists only of finitely many terms and that the free boundary $p=p(t)$ of the price formation problem becomes the zero level set of the heat solution and vice versa. In principle the free boundary can leave and reenter the interval $(-L+a, L-a)$ without impeding the existence of a global solution of \eqref{eq:pfm1}, \eqref{eq:pfm2}, \eqref{eq:initial}, \eqref{eq:newbc}.\\
We shall now construct an 'explicit' solution of \eqref{eq:heat} by separation of variables. We set $F(x,t) = \varphi(x) \psi(t)$ and find $\varphi''/\varphi = \dot{\psi}/\psi = -z^2$. The boundary conditions  give the equation $G(z) := \cos(zL) - \cos(z(L-a)) = 0$ with the eigenfunctions $\varphi(x) = A \sin(zx)$ and $H(z) := \sin(zL)-\sin(z(L-a)) = 0$ with eigenfunctions $\varphi = B \cos(zx)$. We easily compute $G(z) = 0$ iff $z = (2\pi l) / a, (2\pi l)/(2L-a)$ and $H(z) = 0$ iff $z = (2\pi l)/a, (\pi(2l-1)/(2L-a)$ for $l \in \mathbb{Z}$. Next we conclude that $H=H(z)$ is a sine-type function of type $L$ with simple and separated zeros. Therefore $\lbrace \exp(i \frac{2\pi l}{a} x), \exp(i \frac{\pi(2l-1)}{2L-a}x)\rbrace$ is a Riesz basis in $L^2(-L,L)$, see \cite{Avdonin1988}. Also $G = G(z)$ is a sine-type function of type L, its zero are separated except $z=0$, which is a zero of order $2$. Hence $\lbrace \exp(i \frac{2\pi l}{a} x), \exp(i \frac{2 \pi l}{2L-a} x), x\rbrace_{l \in \mathbb{Z}}$ is a Riesz basis in $L^2(-L,L)$, see Theorem D in \cite{Horvath1987}. By separately considering even and odd parts we conclude that $\lbrace\cos(\frac{2\pi l}{a} x), \sin(\frac{2 \pi l}{a} x), \sin(\frac{2\pi l}{2L-a} x), \cos(\frac{\pi (2l-1)}{2L-a} x), 1, x \rbrace$, $l=1,2, \ldots$ is a Riesz basis of eigenfunctions of the heat equation \eqref{eq:heat}. The solution of \eqref{eq:heat} is
\begin{align}\label{eq:solF}
\begin{split}
F(x,t) = \sum_{l=1}^{\infty} &[(A_l \sin(\omega_{1,l} x) + B_l \cos(\omega_{1,l} x)) e^{-\omega_{1,l}^2t} + C_l \sin(\omega_{2,l} x) e^{-\omega_{2,l}^2 t} + D_l \cos(\omega_{3,l}x) e^{-\omega_{3,l}^2t}]{}\\
+ &A_0 x + B_0,
\end{split}
\end{align}
with $\omega_{1,l} = (2\pi l)/a, ~\omega_{2,l} = (2\pi l)/(2L-a),~\omega_{3,l} = ((2l-1)\pi)/(2L-a)$. Note that $A_l, B_l, C_l, D_l, A_0, B_0$ can be determined uniquely for every initial datum $F_I \in L^2(-L,L)$.
\begin{theorem}[Global Existence] The BVP  \eqref{eq:pfm1}, \eqref{eq:pfm2}, \eqref{eq:initial}, \eqref{eq:newbc} has a unique global solution $f = f(x,t)$ for $t>0$. Furthermore the free boundary $p=p(t)$ is a smooth graph $p(t) \in (-L,L)$ for all $t > 0$.
\end{theorem}   
The proof follows from the construction described above. The min-max principle implies that $p=p(t)$ is a graph for $t>0$ and the Hopf principle implies that $p$ is smooth. Finite-time oscillations of $p(t)$ are excluded by the $x$-analyticity of solutions of the heat equations (see \cite{Caffarelli2011}). Since $0 < \int_{-L}^{-L+a} F(x,t) dx$ and $0 > \int_{L-a}^L F(x,t) dx$ are conserved in time (by the equation and the boundary conditions) we conclude that $p(t) \in (-L,L)$ for all $t>0$. \\
Clearly, a solution of \eqref{eq:pfm1}, \eqref{eq:pfm2}, \eqref{eq:initial}, \eqref{eq:newbc} is a solution of \eqref{eq:pfm} on a time interval $[0,T]$ iff $p(t) \in (-L+a, L-a)$ for $t \in [0,T]$. Then the total mass of buyers $M_B = \int_{-\infty}^{p(t)} f(x,t) dx$ and vendors $M_V = -\int_{p(t)}^\infty f(x,t) dx$ are time conserved quantities.
\begin{theorem}\label{t:mass}
The BVP  \eqref{eq:pfm} has a global solution conserving the total mass of buyers and vendors iff the zero level set $p$ of the solution of \eqref{eq:heat} satisfies $p(t) \in (-L+a, L-a)$ for all $t>0$. Then the free boundary $p(t)$ converges to $p_{\infty} \in (-L+a, L-a)$.
\end{theorem}
The proof of Theorem \ref{t:mass} is based on the following lemmas.
\begin{lemma}\label{l:exp}
The solution $F$ converges exponentially fast to $\Fi = \ai{A}x + \ai{B}$ in $L^2(-L,L)$. 
\end{lemma}
{\em Proof:} $~$ From the Riesz base property we deduce that there exist $c_1, c_2 \in \mathbb{R}^+$ such that $c_1 \lVert F_I\rVert^2 \leq \sum_{l=1}^\infty (A_l^2 + B_l^2 + C_l^2 + D_l^2) + A_0^2 + B_0^2  \leq c_2 \lVert F_I \rVert^2$. For $\tilde{F} = F-(A_0x + B_0)$  we deduce
\begin{align*}
c_1 \lVert \tilde{F} \rVert^2 \leq \sum\nolimits_{l=1}^{\infty} [(A_l^2 + B_l^2 + C_l^2+D_l^2)] e^{-2 \gamma_{l} t} \leq e^{-2 \gamma_1 t} \sum\nolimits_{l=1}^{\infty} (A_l^2 + B_l^2 + C_l^2+D_l^2) \leq c_2 e^{-2 \gamma_1 t} \lVert F_I \rVert^2,
\end{align*}
with the sequence $\gamma_l = \min((4\pi^2 l^2)/a^2, (4\pi^2 l^2)/(2L-a)^2, ((2l-1)^2 \pi^2)/(2L-a)^2)$, $l=1,2,\ldots$.\qed
\begin{lemma}\label{l:F}
The solution \eqref{eq:solF} satisfies $\displaystyle\lvert F_x(x,t) \rvert\leq \sup_{x\in(-L,L)} \lvert (F_I)_x \rvert$ and $\lvert F(x,t) \rvert \leq c$, $ \forall x \in (-L,L), ~ t > 0$.
\end{lemma}
{\em Proof:} The function $V=F_x$ satisfies the heat equation with $V(-L,t) = V(-L+a,t),~V(L,t)=V(L-a,t)$ and initial condition $V(x,t=0) = (F_I)_x(x)$. If $V$ assumes its maximum on the cylinder $[-L,L] \times [0,T)$ at either boundary $x = \pm L$, then it must also assume a maximum (with the same value) in the interior (due to bc on $V$). This contradicts the maximum principle, thus $V$ must assume its maximum at $t=0$. The same arguments hold for the minimum. \\
Since $F(x,t) = F(-L,t) + \int_{-L}^x F_x(y,t) dy$ we deduce $a F(-L,t) = -\int_{-L}^{-L+a} F(x,t) dx + \int_{-L}^{-L+a}\int_{-L}^x F_x(y,t) dy dx = I_1 + I_2$.
We know that $F_x$ is bounded, therefore $\lvert I_2 \rvert \leq K$. In addition $\frac{d}{dt} \int_{-L}^{-L+a} F(x,t) dx = 0$ and therefore $\int_{-L}^{-L+a} F(x,t) dx = \int_{-L}^{-L+a} F_I(x) dx$. Thus $F(-L,t)$ as well as $F = F(x,t)$ are bounded uniformly on $(-L,L) \times \mathbb{R}^+$.\qed\\
{\em Proof of Theorem \ref{t:mass}:} We know that $F$ converges exponentially fast to $\Fi = \ai{A} x + \ai{B}$, and that there exists a smooth graph $p = p(t)$ such that $F(p(t),t) = 0$. Now we assume that $p(t) \in (-L+a, L-a)$ for $t>0$, choose any sequence $t_n \rightarrow \infty$ and conclude that there is a subsequence $t_{n_k}$ such that $p(t_{n_k}) \rightarrow p_{\infty} \in [-L+a, L-a]$. Let $\varphi$ be a test function in $\mathcal{D}(p_{\infty},L)$. If $k$ is sufficiently large we conclude $f(x,t_{n_k}) < 0$ for $x \in \supp \varphi$. Therefore $\int_{-L}^L f(x,t_{n_k}) \varphi(x) dx  < 0$. Since $f(\cdot,t) \rightarrow \ai{f}$ there is a subsequence $t_{n_{k_l}}$ such that $f(x,t_{n_{k_l}})$ converges to $\ai{f}(x)$ pointwise a.e. in $(-L,L)$. The function $\lvert f(x,t_{n_k})\rvert \leq K$ on $[-L,L]$ for all $k$. Then we deduce from Lebesgue's' dominated convergence theorem that $\int_{-L}^L f(x,t_{n_k}) \varphi(x) dx \rightarrow \int_{-L}^L f_{\infty} \varphi(x)  \leq 0.$ Since $f_{\infty} = \ai{A} x + \ai{B} + (\ai{A}(x-a) + \ai{B})^- - (\ai{A} (x+a) + \ai{B})^+$ we conclude $\ai{f} \leq 0$ for $x > \ai{p}$ and $\ai{f} \geq 0$ for $x < \ai{p}$. From Lebesgue's' dominated convergence theorem, we conclude $-\int_{p(t_{n_{k_l}})}^L f(x,t) dx \rightarrow -\int_{\ai{p}}^L \ai{f}(x) dx \geq  0.$ Analogously$ \int_{-L}^{\ai{p}} \ai{f} dx \geq 0$. \\
Next we show that $\ai{F}$ has a unique zero in $(-L,L)$. If $\ai{F} = \ai{A} x + \ai{B} \geq 0$ on $(-L,L)$ then $\ai{F}^-=0,~\ai{F}^+=\ai{F}$ and therefore $\ai{f}$ is constant in $(-L,L)$, which is a contradiction. The same argument holds for $\ai{F} \leq 0$. Therefore the function $\ai{f}$ is given by $\ai{f}(x) = \pm \alpha$ for $x \in (-L, \ai{p}-a)$ and $ x\in (\ai{p}+a,L)$ respectively and $\ai{f}= -\alpha/a(x-\ai{p})$ for $x \in [\ai{p}-a, \ai{p}+a]$, with $\alpha \in \mathbb{R}^+$. We conclude that $\ai{p}$ is unique and that $p(t) \rightarrow \ai{p}$ as $t \rightarrow \infty$, since every sequence has a subsequence which converges to the same limit. \qed
\begin{theorem}
Let $f_I$ be such that $M_B/M_V \notin \left[a/(4L-3a), (4L-3a)/a\right]$, where $M_B, M_V$ denotes the initial mass of buyers and vendors. 
Then \eqref{eq:pfm} does not have a global-in-time solution, which conserves both buyers and vendors masses.
\end{theorem}
{\em Proof:} The result follows since $\alpha$, $\ai{p}$ can not be adjusted such that $\ai{p} \in [-L+a,L-a]$, where $M_B = \int_{-L}^{\ai{p}} \ai{f} dx,~M_V = -\int_{\ai{p}}^L \ai{f} dx$. \qed\\
Note that the choice of $L$, which corresponds to the maximally attainable price, is more or less arbitrary but a bad (too small) choice of $L$ might impede global existence. In this case the model clearly looses its 'practical' significance. We remark that global existence results for \eqref{eq:pfm} (with a free boundary which remains in $(-L+a,L-a)$) for initial data which are small perturbations of stationary solutions are straightforward, without using the analytical machinery of \cite{Gonzalez2011}.
\section*{Acknowledgments}
PM acknowledges support by the King Abdullah University of Science and Technology, the Leverhulme Trust and the Royal Society. LC acknowledges support from the Division of Mathematical Sciences of the NSF, MTW  from the Austrian Science Foundation FWF.

\end{document}